\documentclass[12pt]{article}
\usepackage{amsmath,amssymb}

\newtheorem{Theorem}{Theorem}[section]
\newtheorem{Definition}{Definition}[section]
\newtheorem{Corollary}{Corollary}[section]

\newtheorem{Remark}{Remark}

\title{ \Large \bf Coupled coincidence point theorems for mixed $(G,S)$-monotone operators on partially ordered metric spaces and applications}
\author{\small{Habib Yazidi}}
\date{ }
\begin{document}
\maketitle
\abstract
In this paper, we introduce the concept of  mixed $(G,S)$-monotone mappings and prove coupled coincidence and  coupled common fixed point theorems
for such mappings satisfying a nonlinear contraction involving altering distance functions. Presented theorems extend, improve and generalize the very recent results of Harjani, L\'opez and Sadarangani [J. Harjani, B. L\'opez and K. Sadarangani, Fixed point theorems for mixed monotone operators and applications to integral equations, Nonlinear Analysis (2010), doi:10.1016/j.na.2010.10.047] and other existing results in the literature. Some applications to periodic boundary value problems are also considered.

\section{Introduction and preliminaries}
Fixed point problems of contractive mappings in metric spaces endowed with a partially order have been studied by many authors (see \cite{Agarwal}-\cite{Shatanawi}).
Bhaskar and Lakshmikantham \cite{Bhaskar} introduced the concept of a coupled fixed point and studied the problems of a uniqueness of a coupled fixed point in partially ordered metric spaces and applied their theorems to problems of the existence  of solution for a periodic boundary value problem. In \cite{Lakshmikantham}, Lakshmikantham and  \'Ciri\'c established some coincidence and common coupled fixed  point theorems under nonlinear contractions in partially ordered metric spaces.
Very recently, Harjani, L\'opez and Sadarangani \cite{Harjani3} obtained some coupled fixed point theorems for a mixed monotone operator in a complete metric space endowed with a partial order by using altering distance functions. They applied their results to the study of the existence and uniqueness of a nonlinear integral equation.

Now, we briefly recall various basic definitions and facts.

\begin{Definition} (see Bhaskar and Lakshmikantham \cite{Bhaskar}). \label{D1} Let $%
(X,\preceq) $ be a partially ordered set and $F:X\times X\rightarrow X$. Then
the map $F$ is said to have mixed monotone property if $F(x,y)$ is monotone
non-decreasing in $x$ and is monotone non-increasing in $y$, that is, for
any $x,y\in X$,
\begin{equation*}
x_1\preceq x_2\,\,\,\mbox{ implies }\,\,\,F(x_1,y)\preceq F(x_2,y)\,\,\,\mathrm{%
for\,\,\, all}\,\,\,y\in X
\end{equation*}
and
\begin{equation*}
y_1\preceq y_2\,\,\,\mbox{ implies }\,\,\,F(x,y_2)\preceq F(x,y_1)\,\,\,\mathrm{%
for\,\,\, all}\,\,\,x\in X.
\end{equation*}
\end{Definition}

The main result obtained by  Bhaskar and Lakshmikantham \cite{Bhaskar} is the following.
\begin{Theorem}(see Bhaskar and Lakshmikantham \cite{Bhaskar}).
Let $(X,\preceq)$ be a partially ordered set and suppose there is a metric $d$ on $X$ such that $(X,d)$ is a complete metric space.
Let $F: X\times X\rightarrow X$ be a mapping having the mixed monotone property on $X$. Assume that there exists $k\in [0,1)$ such that
$$
d(F(x,y),F(u,v))\leq \frac{k}{2}[d(x,u)+d(y,v)] \mbox{ for each } u\preceq x \mbox{ and } y\preceq v.
$$
Suppose either $F$ is continuous or $X$ has the following properties:
\begin{enumerate}
\item[(i)] if a non-decreasing sequence $x_n\rightarrow x$, then $x_n\preceq x$ for
all $n$,

\item[(ii)] if a non-increasing sequence $x_n\rightarrow x$, then $x\preceq x_n$ for
all $n$.
\end{enumerate}
If there exist $x_0,y_0\in X$ such that
$$
x_0\preceq F(x_0,y_0) \mbox{ and } F(y_0,x_0)\preceq y_0,
$$
then $F$ has a coupled fixed point.
\end{Theorem}

Inspired by Definition \ref{D1}, Lakshmikantham and \'Ciri\'c in \cite{Lakshmikantham}
introduced the concept of a $g$-mixed monotone mapping.

\begin{Definition}(see Lakshmikantham and \'Ciri\'c \cite{Lakshmikantham}). Let $(X,\preceq)$ be a
partially ordered set, $F:X\times X\rightarrow X$ and $g:X\rightarrow X$. Then the map $F$ is
said to have mixed $g$-monotone property if $F(x,y)$ is monotone $g$%
-non-decreasing in $x$ and is monotone $g$-non-increasing in $y$, that is,
for any $x,y\in X$,
\begin{equation*}
gx_1\preceq gx_2\,\,\,\mbox{ implies }\,\,\,F(x_1,y)\preceq F(x_2,y)\,\,\,\mathrm{%
for\,\,\, all}\,\,\,y\in X
\end{equation*}
and
\begin{equation*}
gy_1\preceq gy_2\,\,\,\mbox{ implies }\,\,\,F(x,y_2)\preceq F(x,y_1)\,\,\,\mathrm{%
for\,\,\, all}\,\,\,x\in X.
\end{equation*}
\end{Definition}

\begin{Definition} (see Lakshmikantham and \'Ciri\'c \cite{Lakshmikantham}).
Let $X$ be a non-empty set, and let $F:X\times X\rightarrow X$, $g:X\rightarrow X$ be given mappings.
An element $(x,y)\in X\times X$ is called a coupled coincidence point of the mappings $F$ and $g$ if $F(x,y)=gx$ and  $F(y,x)=gy$.
\end{Definition}

\begin{Definition} (see Lakshmikantham and \'Ciri\'c \cite{Lakshmikantham}). Let $X$ be a
non-empty set. Then we say that the mappings $F:X\times X\rightarrow X$ and $%
g:X\rightarrow X$ are commutative if
\begin{equation*}
g(F(x,y))=F(gx,gy).
\end{equation*}
\end{Definition}
The main result of Lakshmikantham and \'Ciri\'c \cite{Lakshmikantham} is the following.

\begin{Theorem}(see Lakshmikantham and \'Ciri\'c \cite{Lakshmikantham}). Let $(X,\preceq)$ be a
partially ordered set and suppose there is a metric $d$ on $X$ such that $%
(X,d)$ is a complete metric space. Assume there is a function $%
\phi:[0,+\infty)\rightarrow [0,+\infty)$ with $\phi(t)<t$ and $%
\lim_{r\rightarrow t^+}\phi(r)<t$ for each $t>0$ and also suppose $F:X\times
X\rightarrow X$ and $g: X\rightarrow X$ are such that $F$ has the mixed $g$%
-monotone property and
\begin{eqnarray*}
d(F(x,y),F(u,v))\leq\phi\left(\frac{d(gx,gu)+d(gy,gv)}{2}\right)
\end{eqnarray*}
for all $x,y,u,v\in X$ with $gx\preceq gu$ and $gv\preceq gy$. Assume that $%
F(X\times X)\subseteq g(X),$ $g$ is continuous and commutes with $F$ and
also suppose either $F$ is continuous or $X$ has the following properties:

\begin{enumerate}
\item[(i)] if a non-decreasing sequence $x_n\rightarrow x$, then $x_n\preceq x$ for
all $n$,

\item[(ii)] if a non-increasing sequence $x_n\rightarrow x$, then $x\preceq x_n$ for
all $n$.
\end{enumerate}
If there exist $x_0,y_0\in X$ such that $gx_0\preceq F(x_0,y_0)$ and $%
F(y_0,x_0)\preceq gy_0$ then there exist $x,y\in X$ such that $gx=F(x,y)$ and $%
gy=F(y,x)$, that is, $F$ and $g$ have a coupled coincidence point.
\end{Theorem}

Very recently, Harjani, L\'opez and Sadarangani \cite{Harjani3} established coupled fixed point theorems for a mixed monotone operator satisfying contraction involving
altering distance functions in a complete partially ordered metric space. \\
Denote by $\mathcal{F}$ the set of functions $\varphi:[0,+\infty)\rightarrow [0,+\infty)$ satisfying the following properties:\\
(a) $\varphi$ is continuous and non-decreasing,\\
(b) $\varphi(t)=0$ if and only if $t=0$.

\begin{Theorem}\label{THLS}(Harjani, L\'opez and Sadarangani \cite{Harjani3}).
Let $(X,\preceq)$ be a partially ordered set and $d$ be a metric on $X$ such that $(X,d)$ is a complete metric space. Let $F: X\times X\rightarrow X$ be a mapping
having the mixed monotone property on $X$ and satisfying
$$
\varphi(d(F(x,y),F(u,v))\leq \varphi(\max\{d(x,u),d(y,v)\})-\Phi(\max\{d(x,u),d(y,v)\})
$$
for all $x,y,u,v\in X$ with $u\preceq x$ and $y\preceq v$, where $\varphi,\psi \in \mathcal{F}$.
Suppose either $F$ is continuous or $X$ has the following properties:

\begin{enumerate}
\item[(i)] if a non-decreasing sequence $x_n\rightarrow x$, then $x_n\preceq x$ for
all $n$,

\item[(ii)] if a non-increasing sequence $x_n\rightarrow x$, then $x\preceq x_n$ for
all $n$.
\end{enumerate}
If there exist $x_0,y_0\in X$ such that $x_0\preceq F(x_0,y_0)$ and $%
F(y_0,x_0)\preceq y_0$ then $F$ has a coupled fixed point.
\end{Theorem}

In this paper, we introduce the concept of  mixed $(G,S)$-monotone mappings and prove coupled coincidence and  coupled common fixed point theorems
for such mappings satisfying a nonlinear contraction involving altering distance functions. Presented theorems extend, improve and generalize the  results of Harjani, L\'opez and Sadarangani \cite{Harjani3}. As applications of our obtained results, we study the existence and uniqueness of  solution to periodic boundary value problem.

\section{Main Results}
Now, we introduce the concept of mixed $(G,S)$-monotone property.

\begin{Definition}
Let $X$ be a non-empty set endowed with a partial order $\preceq$. Consider the mappings $F: X\times X\rightarrow X$ and $G,S: X\rightarrow X$.
We say that $F$ has the mixed $(G,S)$-monotone property on $X$ if for all $x,\,y\,\in X$,
\begin{eqnarray*}
x_{1},\,x_{2}\,\in X,\quad G(x_{1})\preceq S(x_{2})\Rightarrow F(x_{1},y)\preceq F(x_{2},y),\\
x_{1},\,x_{2}\,\in X,\quad G(x_{1})\succeq S(x_{2})\Rightarrow F(x_{1},y)\succeq F(x_{2},y),\\
y_{1},\,y_{2}\,\in X,\quad G(y_{1})\preceq S(y_{2})\Rightarrow F(x,y_{1})\succeq F(x,y_{2}),\\
y_{1},\,y_{2}\,\in X,\quad G(y_{1})\succeq S(y_{2})\Rightarrow F(x,y_{1})\preceq F(x,y_{2}).
\end{eqnarray*}
\end{Definition}
\begin{Remark}
If we take $G=S$, then $F$ has the mixed $(G,S)$-monotone property implies that $F$ has the mixed $G$-monotone property.
\end{Remark}

Now, we state and prove our first result.

\begin{Theorem}\label{T1}
Let $(X,\,\preceq)$ be a partially ordered set and suppose that there exists a metric $d$ on $X$ such that $(X,d)$ is a complete metric space. Let $G,S: X\rightarrow X$ and $F: X\times X\rightarrow X$ be a mapping having the mixed $(G,S)$-monotone property on $X$. Suppose that
\begin{equation}
\varphi(d(F(x,y),F(u,v)))\leq \varphi(\max\{d(Gx,Su),d(Sy,Gv)\})-\phi(\max \{d(Gx,Su),d(Sy,Gv)\}),
\label{eq1}
\end{equation}
for all $x,\,y,\,u,\,v\,\in X$ with $G(x)\preceq S(u)$ or $G(x)\succeq S(u)$ and $S(y)\succeq G(v)$ or $S(y)\preceq G(v)$, where $\varphi, \phi\in \mathcal{F}$.
Assume that $F(X\times X)\subseteq G(X)\cap S(X)$ and assume also that $G,S$ and $F$ satisfy the following hypotheses:
\begin{enumerate}
\item[(I)] $F, G$ and $S$ are continuous,
\item[(II)] $F$ commutes respectively with $G$ and $S$.
\end{enumerate}
If there exist $x_{0},\,y_{0},\,x_{1}$ and $y_{1}$ such that
\begin{equation*}
\left\{\begin{array}{ll}
G(x_{0})\preceq S(x_{1})\preceq F(x_{0},y_{0});\\
G(y_{0})\succeq S(y_{1})\succeq F(y_{0},x_{0}),
\end{array}
\right.
\end{equation*}
then there exist $x,\,y\in X$ such that
$$G(x)=S(x)=F(x,y)\quad\textrm{and}\quad G(y)=S(y)=F(y,x) ,$$
 that is, $G,S$ and $F$ have a coupled coincidence point $(x,y)\in X\times X$.
\end{Theorem}

\noindent {\bf Proof.} Let $x_{0},\,y_{0},\,x_{1},\,y_{1}\in\, X$ such that
$$
G(x_{0})\preceq S(x_{1})\preceq F(x_{0},y_{0}) \quad \mbox{ and }\quad  G(y_{0})\succeq S(y_{1})\succeq F(y_{0},x_{0}).
$$
Since $F(X\times X)\subseteq G(X)\cap S(X)$, we can choose $x_{2},\,y_{2},\,x_{3},\,y_{3}\,\in X$ such that
\begin{equation*}
\begin{array}{ccc}
\left\{\begin{array}{ll}
G(x_{2})=F(x_{0},y_{0})\\
G(y_{2})=F(y_{0},x_{0})
\end{array}
\right.
\quad\textrm{and}\quad
\left\{\begin{array}{ll}
S(x_{3})=F(x_{1},y_{1})\\
S(y_{3})=F(y_{1},x_{1})
\end{array}
\right.\cdot
\end{array}
\end{equation*}
Continuing this process we can construct sequences $\{x_{n}\}$ and $\{y_{n}\}$ in $X$ such that
\begin{equation}
\begin{array}{ccc}
\left\{\begin{array}{ll}
G(x_{2n+2})=F(x_{2n},y_{2n})\\
G(y_{2n+2})=F(y_{2n},x_{2n})
\end{array}
\right.;\,\,
\left\{\begin{array}{ll}
S(x_{2n+3})=F(x_{2n+1},y_{2n+1})\\
S(y_{2n+3})=F(y_{2n+1},x_{2n+1})
\end{array}
\right.
\quad\textrm{for all $n\geq 0$.}
\end{array}
\end{equation}

We shall show that for all $n\geq 0$,
\begin{equation}\label{eq2}
G(x_{2n})\preceq S(x_{2n+1})\preceq G(x_{2n+2})
\end{equation}
and
\begin{equation}\label{eq222}
G(y_{2n})\geq S(y_{2n+1})\geq G(y_{2n+2}).
\end{equation}
As
$G(x_{0})\preceq S(x_{1})\preceq F(x_{0},y_{0})=G(x_{2})$ and $G(y_{0})\succeq S(y_{1})\succeq F(y_{0},x_{0})=G(y_{2})$,  our claim is satisfied for $n=0$.\\
Suppose that (\ref{eq2}) and (\ref{eq222}) hold for some fixed $n\geq 0$. Since $G(x_{2n})\preceq S(x_{2n+1})\preceq G(x_{2n+2})$ and $G(y_{2n})\succeq S(y_{2n+1})\succeq G(y_{2n+2})$, and as $F$ has the mixed $(G,S)$-monotone property, we have
$$
G(x_{2n+2})=F(x_{2n},y_{2n})\preceq F(x_{2n+1},y_{2n})\preceq F(x_{2n+1},y_{2n+1})\preceq F(x_{2n+2},y_{2n+1})\preceq F(x_{2n+2},y_{2n+2}),
$$
then
$$
G(x_{2n+2})\preceq S(x_{2n+3})\preceq G(x_{2n+4}).
$$
On the other hand,
$$
G(y_{2n+2})=F(y_{2n},x_{2n})\succeq F(y_{2n+1},x_{2n})\succeq F(y_{2n+1},x_{2n+1})\succeq F(y_{2n+2},x_{2n+1})\succeq F(y_{2n+2},x_{2n+2}),
$$
then
$$
G(y_{2n+2})\succeq S(y_{2n+3})\succeq G(y_{2n+4}).
$$
Thus by induction, we proved that (\ref{eq2}) and (\ref{eq222}) hold for all $n\geq 0$.

We complete the proof in the following steps\\

\noindent {\bf Step 1:} We will prove that
\begin{equation}
\lim_{n\rightarrow +\infty}d(F(x_{n},y_{n}),\,F(x_{n+1},y_{n+1}))=\lim_{n\rightarrow +\infty}d(F(y_{n},x_{n}),\,F(y_{n+1},x_{n+1}))=0 .
\label{eq4}
\end{equation}
From (\ref{eq2}), (\ref{eq222}) and (\ref{eq1}), we have
\begin{align}
&\varphi(d(F(x_{2n},y_{2n}),F(x_{2n+1},y_{2n+1})))\\
&\leq \varphi(\max \{d(G x_{2n},S x_{2n+1}),d(G y_{2n},S y_{2n+1})\})-\phi(\max \{d(G x_{2n},S x_{2n+1}),d(G y_{2n},S y_{2n+1})\})\notag\\
&\leq \varphi(\max \{d(G x_{2n},S x_{2n+1}),d(G y_{2n},S y_{2n+1})\}).
\label{bla1}
\end{align}
Since $\varphi$ is a non-decreasing function, we get that
\begin{equation*}
d(F(x_{2n},y_{2n}),F(x_{2n+1},y_{2n+1}))\leq \max \{d(G x_{2n},S x_{2n+1}),d(G y_{2n},S y_{2n+1})\}.
\end{equation*}
Therefore
\begin{equation}
d(Gx_{2n+2},Sx_{2n+3})\leq \max \{d(G x_{2n},S x_{2n+1}),d(G y_{2n},S y_{2n+1})\}.
\label{eq5}
\end{equation}
Again,  using (\ref{eq2}), (\ref{eq222}) and (\ref{eq1}), we have
\begin{align}
\notag &\varphi(d(F(y_{2n},x_{2n}),F(y_{2n+1},x_{2n+1})))\\
\notag &\leq \varphi(\max \{d(G y_{2n},S y_{2n+1}),d(G x_{2n},S x_{2n+1})\})-\phi(\max \{d(G y_{2n},S y_{2n+1}),d(G x_{2n},S x_{2n+1})\})\\
&\leq \varphi(\max \{d(G y_{2n},S y_{2n+1}),d(G x_{2n},S x_{2n+1})\}).
\label{bla2}
\end{align}
Since $\varphi$ is non-decreasing, we have
\begin{equation*}
d(F(y_{2n},x_{2n}),F(y_{2n+1},x_{2n+1}))\leq \max \{d(G y_{2n},S y_{2n+1}),d(G x_{2n},S x_{2n+1}\}.
\end{equation*}
Therefore
\begin{equation}
d(Gy_{2n+2},Sy_{2n+3})\leq \max \{d(G y_{2n},S y_{2n+1}),d(G x_{2n},S x_{2n+1})\}.
\label{eq6}
\end{equation}
Combining (\ref{eq5}) and (\ref{eq6}), we obtain
\begin{equation*}
\max\{d(Gx_{2n+2},Sx_{2n+3}),d(Gy_{2n+2},Sy_{2n+3})\}\leq \max \{d(G x_{2n},S x_{2n+1}),d(G y_{2n},S y_{2n+1})\}.
\end{equation*}
Then $\bigg\{\max \{d(G x_{2n},S x_{2n+1}),d(G y_{2n},S y_{2n+1})\}\bigg\}$ is a positive decreasing sequence. Hence there exists
$r\geq 0$ such that
$$
\lim_{n\rightarrow +\infty}\max \{d(G x_{2n},S x_{2n+1}),d(G y_{2n},S y_{2n+1})\}=r.
$$
Combining (\ref{bla1}) and (\ref{bla2}), we obtain
\begin{align*}
&\max\{\varphi(d(G x_{2n+2},S x_{2n+3})),\varphi(d(G y_{2n+2},S y_{2n+3}))\}\\
&\leq \varphi(\max \{d(G x_{2n},S x_{2n+1}),d(G y_{2n},S y_{2n+1})\})-\phi(\max \{d(G x_{2n},S x_{2n+1}),d(G y_{2n},S y_{2n+1})\}).
\end{align*}
Since $\varphi$ is non-decreasing, we get
\begin{align*}
&\varphi(\max\{d(G x_{2n+2},S x_{2n+3}),d(G y_{2n+2},S y_{2n+3})\})\\
&\leq \varphi(\max \{d(G x_{2n},S x_{2n+1}),d(G y_{2n},S y_{2n+1})\})-\phi(\max \{d(G x_{2n},S x_{2n+1}),d(G y_{2n},S y_{2n+1})\}).
\end{align*}
Letting $n\rightarrow +\infty$ in the above inequality, we get
$$
\varphi(r)\leq \varphi(r)-\phi(r),
$$
which implies that $\phi(r)=0$ and then, since $\phi$ is an altering distance function, $r=0$.
Consequently
\begin{equation}\label{bs1}
\lim_{n\rightarrow +\infty }\max \{d(F(x_{2n},y_{2n}),F(x_{2n+1},y_{2n+1})),d(F(y_{2n},x_{2n}),F(y_{2n+1},x_{2n+1}))\}=0.
\end{equation}
By the same way, we obtain
\begin{equation}\label{bs2}
\lim_{n\rightarrow +\infty }\max \{d(F(x_{2n+1},y_{2n+1}),F(x_{2n+2},y_{2n+2})),d(F(y_{2n+1},x_{2n+1}),F(y_{2n+2},x_{2n+2}))\}=0.
\end{equation}
Finally, (\ref{bs1}) and (\ref{bs2}) give the desired  result, that is, (\ref{eq4}) holds.\\

\noindent {\bf Step 2:} We will prove that $F(x_{n},y_{n})$ and $F(y_{n},x_{n})$ are Cauchy sequences.\\
From (\ref{eq4}), it is sufficient to show that $F(x_{2n},y_{2n})$ and $F(y_{2n},x_{2n})$ are Cauchy sequences.\\
We proceed by negation and suppose that at least one of the sequences $F(x_{2n},y_{2n})$ or $F(y_{2n},x_{2n})$ is not a Cauchy sequence.\\
This implies that $d(F(x_{2n},y_{2n}),F(x_{2m},y_{2m}))\nrightarrow 0$ or $d(F(y_{2n},x_{2n}),F(y_{2m},x_{2m}))\nrightarrow 0$ as $n,m\rightarrow +\infty$.\\
Consequently
$$
\max\{d(F(x_{2n},y_{2n}),F(x_{2m},y_{2m})),d(F(y_{2n},x_{2n}),F(y_{2m},x_{2m}))\}\nrightarrow 0,\,\,\,\textrm{ as }  n,m\rightarrow +\infty.
$$
Then there exists $\varepsilon>0$ for which we can find two subsequences of positive integers $\{m(i)\}$ and $\{n(i)\}$ such that $n(i)$ is the smallest index for which $n(i)>m(i)>i$,
\begin{equation}\label{C1}
\max\{d(F(x_{2m(i)},y_{2m(i)}),F(x_{2n(i)},y_{2n(i)})),d(F(y_{2m(i)},x_{2m(i)}),F(y_{2n(i)},x_{2n(i)}))\}\geq\varepsilon.
\end{equation}
This means that
\begin{equation}\label{C2}
\max\{d(F(x_{2m(i)},y_{2m(i)}),F(x_{2n(i)-2},y_{2n(i)-2})),d(F(y_{2m(i)},x_{2m(i)}),F(y_{2n(i)-2},x_{2n(i)-2}))\}<\varepsilon.
\end{equation}
From (\ref{C1}), (\ref{C2}) and using the triangular inequality, we get
\begin{align*}
&\varepsilon \leq \max\{d(F(x_{2m(i)},y_{2m(i)}),F(x_{2n(i)},y_{2n(i)})),d(F(y_{2m(i)},x_{2m(i)}),F(y_{2n(i)},x_{2n(i)}))\}\\
&\leq \max\{d(F(x_{2m(i)},y_{2m(i)}),F(x_{2n(i)-2},y_{2n(i)-2})),d(F(y_{2m(i)},x_{2m(i)}),F(y_{2n(i)-2},x_{2n(i)-2}))\}\\
&+\max\{d(F(x_{2n(i)-2},y_{2n(i)-2}),F(x_{2n(i)-1},y_{2n(i)-1})),d(F(y_{2n(i)-2},x_{2n(i)-2}),F(y_{2n(i)-1},x_{2n(i)-1}))\}\\
&+\max\{d(F(x_{2n(i)-1},y_{2n(i)-1}),F(x_{2n(i)},y_{2n(i)})),(F(y_{2n(i)-1},x_{2n(i)-1}),F(y_{2n(i)},x_{2n(i)}))\}\\
&<\varepsilon+\max\{d(F(x_{2n(i)-2},y_{2n(i)-2}),F(y_{2n(i)-1},x_{2n(i)-1})),d(F(y_{2n(i)-2},x_{2n(i)-2}),F(y_{2n(i)-1},x_{2n(i)-1}))\}\\
&+\max\{(F(x_{2n(i)-1},y_{2n(i)-1}),F(x_{2n(i)},y_{2n(i)})),d(F(y_{2n(i)-1},x_{2n(i)-1}),F(y_{2n(i)},x_{2n(i)}))\}.
\end{align*}
Letting $i\rightarrow+\infty$ in above inequality and using (\ref{eq4}), we obtain that
\begin{equation}\label{L1}
\lim_{i\rightarrow+\infty}\max(d(F(x_{2m(i)},y_{2m(i)}),F(x_{2n(i)},y_{2n(i)})),d(F(y_{2m(i)},x_{2m(i)}),F(y_{2n(i)},x_{2n(i)})))=\varepsilon.
\end{equation}
Also, we have
\begin{align*}
&\varepsilon \leq \max\{d(F(x_{2m(i)},y_{2m(i)}),F(x_{2n(i)},y_{2n(i)})),d(F(y_{2m(i)},x_{2m(i)}),F(y_{2n(i)},x_{2n(i)}))\}\\
&\leq \max\{d(F(x_{2m(i)},y_{2m(i)}),F(y_{2m(i)-1},x_{2m(i)-1})),d(F(y_{2m(i)},x_{2m(i)}),F(y_{2m(i)-1},x_{2m(i)-1}))\}\\
&+\max\{d(F(x_{2m(i)-1},y_{2m(i)-1}),F(x_{2n(i)},y_{2n(i)})),d(F(y_{2m(i)-1},x_{2m(i)-1}),F(y_{2n(i)},x_{2n(i)}))\}\\
&\leq 2\max\{d(F(x_{2m(i)},y_{2m(i)}),F(x_{2m(i)-1},y_{2m(i)-1})),d(F(y_{2m(i)},x_{2m(i)}),F(y_{2m(i)-1},x_{2m(i)-1}))\}\\
&+\max\{d(F(x_{2m(i)},y_{2m(i)}),F(x_{2n(i)},y_{2n(i)})),d(F(y_{2m(i)},x_{2m(i)}),F(y_{2n(i)},x_{2n(i)}))\}.
\end{align*}
Using (\ref{eq4}), (\ref{L1}) and letting $i\rightarrow +\infty$ in the above inequality, we obtain
\begin{equation}\label{L2}
\lim_{i\rightarrow+\infty}\max\{d(F(x_{2m(i)-1},y_{2m(i)-1}),F(x_{2n(i)},y_{2n(i)})),d(F(y_{2m(i)-1},x_{2m(i)-1}),F(y_{2n(i)},x_{2n(i)}))\}=\varepsilon.
\end{equation}
On other hand, we have
\begin{align*}
&\max\{d(F(x_{2m(i)},y_{2m(i)}),F(x_{2n(i)},y_{2n(i)})),d(F(y_{2m(i)},x_{2m(i)}),F(y_{2n(i)},x_{2n(i)}))\}\\
&\leq \max\{d(F(x_{2m(i)},y_{2m(i)}),F(x_{2n(i)+1},y_{2n(i)+1})),d(F(y_{2m(i)},x_{2m(i)}),F(y_{2n(i)+1},x_{2n(i)+1}))\}\\
&+\max\{d(F(x_{2n(i)+1},y_{2n(i)+1}),F(x_{2n(i)},y_{2n(i)})),d(F(y_{2n(i)+1},x_{2n(i)+1}),F(y_{2n(i)},x_{2n(i)}))\}.
\end{align*}
Since $\varphi$ is a continuous non-decreasing function, it follows from the above inequality that
\begin{align}\label{eq7}
&\varphi(\varepsilon)\leq \\
\notag & \limsup_{i\rightarrow+\infty}\varphi(\max\{d(F(x_{2m(i)},y_{2m(i)}),F(x_{2n(i)+1},y_{2n(i)+1})),d(F(y_{2m(i)},x_{2m(i)}),F(y_{2n(i)+1},x_{2n(i)+1}))\}).
\end{align}
Using the contractive condition, on one hand we have
\begin{align*}
&\varphi(d(F(x_{2m(i)},y_{2m(i)})),F(x_{2n(i)+1},y_{2n(i)+1}))\leq \varphi(\max\{d(G x_{2m(i)},S x_{2n(i)+1}),d(G y_{2m(i)},S y_{2n(i)+1})\})\\
&-\phi(\max\{d(G x_{2m(i)},S x_{2n(i)+1}),d(G y_{2m(i)},Sy_{2n(i)+1})\}) \leq\\ &\varphi(\max\{d(F(x_{2m(i)-2},y_{2m(i)-2})),F(x_{2n(i)-1},y_{2n(i)-1}),d(F(y_{2m(i)-2},x_{2m(i)-2})),F(y_{2n(i)-1},x_{2n(i)-1})\})\\
&-\phi(\max\{d(F(x_{2m(i)-2},y_{2m(i)-2})),F(x_{2n(i)-1},y_{2n(i)-1}),d(F(y_{2m(i)-2},x_{2m(i)-2})),F(y_{2n(i)-1},x_{2n(i)-1})\}).
\end{align*}
On the other hand we have
\begin{align*}
&\varphi(d(F(y_{2m(i)},x_{2m(i)})),F(y_{2n(i)+1},x_{2n(i)+1}))\leq \varphi(\max\{d(G y_{2m(i)},S y_{2n(i)+1}),d(G x_{2m(i)},S x_{2n(i)+1})\})\\
&-\phi(\max\{d(G y_{2m(i)},S y_{2n(i)+1}),d(G x_{2m(i)},S x_{2n(i)+1})\})\leq\\
&\varphi(\max\{d(F(y_{2m(i)-2},x_{2m(i)-2})),F(y_{2n(i)-1},x_{2n(i)-1}),d(F(x_{2m(i)-2},y_{2m(i)-2})),F(x_{2n(i)-1},y_{2n(i)-1}))\})\\
&-\phi(\max\{d(F(y_{2m(i)-2},x_{2m(i)-2})),F(y_{2n(i)-1},x_{2n(i)-1}),d(F(x_{2m(i)-2},y_{2m(i)-2})),F(x_{2n(i)-1},y_{2n(i)-1})\}).
\end{align*}
Therefore
\begin{align}\label{eq8}
\notag &\max\{\varphi(d(F(x_{2m(i)},y_{2m(i)})),F(x_{2n(i)+1},y_{2n(i)+1})),\varphi(d(F(y_{2m(i)},x_{2m(i)})),F(y_{2n(i)+1},x_{2n(i)+1}))\}\\
\notag &\leq\varphi(\max\{d(G x_{2m(i)},S x_{2n(i)+1}),d(G y_{2m(i)},S y_{2n(i)+1})\})
\notag -\phi(\max\{d(F(x_{2m(i)-2},y_{2m(i)-2})),\\
&\quad F(x_{2n(i)-1},y_{2n(i)-1}),d(F(y_{2m(i)-2},x_{2m(i)-2})),F(y_{2n(i)-1},x_{2n(i)-1})\}).
\end{align}
We claim that
\begin{align}\label{eq9}
\notag &\max\{d(F(x_{2m(i)-2},y_{2m(i)-2})),F(x_{2n(i)-1},y_{2n(i)-1}),d(F(y_{2m(i)-2},x_{2m(i)-2})),F(y_{2n(i)-1},x_{2n(i)-1})\})\\
&\rightarrow\varepsilon \,\,\textrm{ as }\,\, i\rightarrow +\infty.
\end{align}
In fact, using the triangular inequality, we have
\begin{align*}
&d(F(x_{2m(i)-2},y_{2m(i)-2}),F(x_{2n(i)-1},y_{2n(i)-1}))\\
&\leq d(F(x_{2m(i)-2},y_{2m(i)-2}),F(x_{2m(i)-1},y_{2m(i)-1}))+d(F(x_{2m(i)-1},y_{2m(i)-1}),F(x_{2n(i)},y_{2n(i)}))\\
&+d(F(x_{2n(i)},y_{2n(i)}),F(x_{2n(i)-1},y_{2n(i)-1})).
\end{align*}
Letting $i\rightarrow +\infty$ in the above inequality and using (\ref{eq4}) and (\ref{L2}), we obtain
\begin{equation}\label{eq10}
\lim_{i\rightarrow +\infty}d(F(x_{2m(i)-2},y_{2m(i)-2}),F(x_{2n(i)-1},y_{2n(i)-1}))\leq\varepsilon.
\end{equation}
On the other hand, we have
\begin{align*}
&d(F(x_{2m(i)-1},y_{2m(i)-1}),F(x_{2n(i)},y_{2n(i)}))\\
&\leq d(F(x_{2m(i)-1},y_{2m(i)-1}),F(x_{2m(i)-2},y_{2m(i)-2}))+d(F(x_{2m(i)-2},y_{2m(i)-2}),F(x_{2n(i)-1},y_{2n(i)-1}))\\
&+d(F(x_{2n(i)-1},y_{2n(i)-1}),F(x_{2n(i)},y_{2n(i)})).
\end{align*}
Letting $i\rightarrow +\infty$ in the above inequality and using (\ref{eq4}) and (\ref{L2}), we obtain
\begin{equation}\label{eq11}
\varepsilon\leq\lim_{i\rightarrow +\infty}d(F(x_{2m(i)-2},y_{2m(i)-2}),F(x_{2n(i)-1},y_{2n(i)-1})).
\end{equation}
Combining (\ref{eq10}) and (\ref{eq11}), we get
$$
\lim_{i\rightarrow +\infty}d(F(x_{2m(i)-2},y_{2m(i)-2}),F(x_{2n(i)-1},y_{2n(i)-1}))=\varepsilon.
$$
By the same way, we obtain
$$
\lim_{i\rightarrow +\infty}d(F(y_{2m(i)-2},x_{2m(i)-2}),F(y_{2n(i)-1},x_{2n(i)-1}))=\varepsilon.
$$
Thus we proved (\ref{eq9}). Finally, letting $i\rightarrow +\infty$ in (\ref{eq8}), using (\ref{eq7}), (\ref{eq9}) and the continuity of $\varphi$ and $\phi$, we get
$\varphi(\varepsilon)\leq \varphi(\varepsilon)-\phi(\varepsilon)$,
which implies that $\phi(\varepsilon)=0$, that is,  $\varepsilon=0$, a contradiction.
Thus $(F(x_{2n},y_{2n}))$ and $(F(y_{2n},x_{2n}))$ are Cauchy sequences in $X$, which gives us that $(F(x_{n},y_{n}))$ and $(F(y_{n},x_{n}))$ are also Cauchy sequences.\\

\noindent {\bf Step 3:} Existence of a coupled coincidence point.\\
Since$(F(x_{n},y_{n}))$ and $((F(y_{n},x_{n})))$ are Cauchy sequences in the complete metric space $(X,\,d)$, there exist $\alpha,\,\alpha'\in X$ such that:
$$\lim_{n\rightarrow +\infty}F(x_{n},y_{n})=\alpha\quad\textrm{and}\quad \lim_{n\rightarrow +\infty}F(y_{n},x_{n})=\alpha'.$$
Therefore,  $\displaystyle \lim_{n\rightarrow +\infty}G(x_{2n+2})=\alpha$, $\displaystyle \lim_{n\rightarrow +\infty}G(y_{2n+2})=\alpha'$, $\displaystyle \lim_{n\rightarrow +\infty}S(x_{2n+3})=\alpha$ and $\displaystyle \lim_{n\rightarrow +\infty}S(y_{2n+3})=\alpha'$.\\
using the continuity and the commutativity of $F$ and $G$, we have
\begin{equation*}
\begin{array}{ccc}
\begin{array}{lll}
G(G(x_{2n+2}))&=&G(F(x_{2n},y_{2n}))\\
&=&F(Gx_{2n},Gy_{2n})
\end{array}
\quad\textrm{and}\quad
\begin{array}{lll}
G(G(y_{2n+2}))&=&G(F(y_{2n},x_{2n}))\\
&=&F(Gy_{2n},Gx_{2n}).
\end{array}
\end{array}
\end{equation*}
Letting $n\rightarrow +\infty$, we get
$G(\alpha)=F(\alpha,\alpha')$ and $G(\alpha')=F(\alpha',\alpha).$\\
Using also the continuity and the commutativity of $F$ and $S$, by the same way, we obtain $S(\alpha)=F(\alpha,\alpha')$ and $S(\alpha')=F(\alpha',\alpha)$.\\
Therefore
$$G(\alpha)=F(\alpha,\alpha')=S(\alpha)\quad\textrm{and}\quad G(\alpha')=F(\alpha',\alpha)=S(\alpha').$$
Thus we proved that $(\alpha,\alpha')$ is a coupled coincidence point of $G, S$ and $F$.
\hfill  $\blacksquare$

In the next result, we prove that the previous theorem is still valid if we replace the continuity of $F$ by some conditions.

\begin{Theorem}\label{T2}
If we replace the continuity hypothesis of $F$ in Theorem \ref{T1} by the following conditions:
\begin{itemize}
\item[(i)] if $(x_{n})$ is a non-decreasing sequences with $x_{n}\rightarrow x$ then $x_{n}\leq x$ for each $n\in \mathbb{N}$,
\item[(ii)] if $(y_{n})$ is a non-increasing sequences with $y_{n}\rightarrow y$ then $y\leq y_{n}$ for each $n\in \mathbb{N}$,
\item[(iii)] $x,y\in X,\quad x\preceq y \Rightarrow Gx\preceq Sy$,
\item[(iv)] $x,y\in X,\quad x\succeq y \Rightarrow Gx\succeq Sy$.
\end{itemize}
Then  $G,S$ and $F$ have a coupled coincidence point.
\end{Theorem}

\noindent {\bf Proof.}
Following the proof of Theorem \ref{T1}, we have that $F(x_{n},y_{n})$ and $F(y_{n},x_{n})$ are Cauchy sequences in the complete metric space $(X,d)$, there exist $\alpha$, $\alpha'\in X$ such that
$$
\lim_{n\rightarrow +\infty}F(x_{n},y_{n})=\alpha\quad\textrm{and}\quad \lim_{n\rightarrow +\infty}F(y_{n},x_{n})=\alpha'.
$$
Therefore $ \displaystyle\lim_{n\rightarrow +\infty}F(x_{2n},y_{2n})=\alpha$ and $\displaystyle\lim_{n\rightarrow +\infty}F(y_{2n},x_{2n})=\alpha'.$ Hence $\displaystyle \lim_{n\rightarrow +\infty}G(x_{2n+2})=\alpha$, $\displaystyle \lim_{n\rightarrow +\infty}G(y_{2n+2})=\alpha'$, $\displaystyle \lim_{n\rightarrow +\infty}S(x_{2n+3})=\alpha$ and $\displaystyle \lim_{n\rightarrow +\infty}S(y_{2n+3})=\alpha'$.
Using the commutativity of $F$ and $G$ and of $F$ and $S$ and the contractive condition, it follows from conditions (iii)-(iv) that
\begin{align}\label{eq12}
\notag&\varphi(d(G(F(x_{2n},y_{2n})),S(F(x_{2n+1},y_{2n+1}))))\\
\notag &=\varphi(d(F(G x_{2n},G y_{2n}),F(S x_{2n+1}, S y_{2n+1})))\\
&\leq \varphi(\max\{d(G(G x_{2n}),S(S x_{2n+1})),d(G(G y_{2n}),S(Sy_{2n+1}))\})\\
&-\phi(\max\{d(G(Gx_{2n}),S(Sx_{2n+1})),d(G(Gy_{2n}),S(Sy_{2n+1}))\}).
\end{align}
Similarly, we have
\begin{align}\label{eq13}
\notag&\varphi(d(G(F(y_{2n},x_{2n})),S(F(y_{2n+1},x_{2n+1}))))\\
\notag &=\varphi(d(F(G y_{2n},G x_{2n}),F(S y_{2n+1}, S x_{2n+1})))\\
&\leq \varphi(\max\{d(G(Gy_{2n}),S(Sy_{2n+1})),d(G(G x_{2n}),S(Sx_{2n+1}))\})\\
&-\phi(\max\{d(G(Gy_{2n}),S(Sy_{2n+1})),d(G(Gx_{2n}),S(Sx_{2n+1}))\}).
\end{align}
Combining (\ref{eq12}), (\ref{eq13}) and the fact that $\max\{\varphi(a),\varphi(b)\}=\varphi(\max\{a,b\})$ for $a,b\in[0,+\infty)$, from (iii)-(iv), we obtain
\begin{align*}
&\varphi(\max\{d(G(F(x_{2n},y_{2n})),S(F(x_{2n+1},y_{2n+1}))),d(G(F(y_{2n},x_{2n})),S(F(y_{2n+1},x_{2n+1})))\})\\
&\leq \varphi(\max\{d(G(G x_{2n}),S(Sx_{2n+1})),d(G(Gy_{2n}),S(Sy_{2n+1}))\})\\
&-\phi(\max\{d(G(Gx_{2n}),S(Sx_{2n+1})),d(G(Gy_{2n}),S(Sy_{2n+1}))\}).
\end{align*}
Letting $n\rightarrow +\infty$ in the last expression, using the continuity of $G$ and $S$, we get
\begin{align*}
&\varphi(\max\{d(G(\alpha),S(\alpha)),d(G(\alpha'),S(\alpha'))\})\\
&\leq \varphi(\max\{d(G(\alpha),S(\alpha)),d(G(\alpha'),S(\alpha'))\})-\phi(\max\{d(G(\alpha),S(\alpha)),d(G(\alpha'),S(\alpha'))\}).
\end{align*}
This implies that $\phi(\max\{d(G(\alpha),S(\alpha)),d(G(\alpha'),S(\alpha'))\})=0$ and, since $\phi$ is an altering distance function, then
\begin{equation*}
\max\{d(G(\alpha),S(\alpha)),d(G(\alpha'),S(\alpha'))\}=0.
\end{equation*}
Consequently
\begin{eqnarray}\label{eq14}
G(\alpha)=S(\alpha) \quad\textrm{and}\quad G(\alpha')=S(\alpha').
\end{eqnarray}
To finish the proof, we claim that $F(\alpha,\alpha')=G(\alpha)=S(\alpha)$ and $F(\alpha',\alpha)=G(\alpha')=S(\alpha')$.\\
Indeed, using the contractive condition, it follows from (i)-(iv) that
\begin{align*}
&\varphi(d(F(G x_{2n},G y_{2n}),F(\alpha,\alpha')))\\
&\leq \varphi(\max\{d(G(G x_{2n}),S(\alpha)),d(G(G y_{2n}),S(\alpha'))\})-\phi(\max\{d(G(G x_{2n}),S(\alpha)),d(G(G y_{2n}),S(\alpha'))\})\\
&\leq \varphi(\max\{d(G(G x_{2n}),S(\alpha)),d(G(G y_{2n}),S(\alpha'))\}).
\end{align*}
Using the fact that $\varphi$ is non-decreasing, we get
\begin{equation}\label{eq15}
d(F(G x_{2n},G y_{2n}),F(\alpha,\alpha'))\leq \max\{d(G(G x_{2n}),S(\alpha)),d(G(G y_{2n}),S(\alpha'))\}.
\end{equation}
Similarly, we have
\begin{eqnarray*}
\varphi(d(F(G y_{2n},G x_{2n}),F(\alpha',\alpha)))&\leq \varphi(\max\{d(G(G y_{2n}),S(\alpha')),d(G(G x_{2n}),S(\alpha))\})\\
&-\phi(\max\{d(G(G y_{2n}),S(\alpha')),d(G(G x_{2n}),S(\alpha))\}\\
&\leq\varphi(\max\{d(G(G y_{2n}),S(\alpha')),d(G(G x_{2n}),S(\alpha))\}).
\end{eqnarray*}
Using the fact that $\varphi$ is non-decreasing, we see that
\begin{equation}\label{eq16}
d(F(G y_{2n},G x_{2n}),F(\alpha',\alpha))\leq \max\{d(G(G y_{2n}),S(\alpha')),d(G(G x_{2n}),S(\alpha))\}.
\end{equation}
Combining (\ref{eq15}) and (\ref{eq16}), we get
\begin{align*}
&\max\{d(F(G x_{2n},G y_{2n}),F(\alpha,\alpha')),d(F(G y_{2n},G x_{2n}),F(\alpha',\alpha)))\\
&\leq \max\{d(G(G y_{2n}),S(\alpha')),d(G(G x_{2n}),S(\alpha))\}.
\end{align*}
Using the commutativity of $F$ and $G$, we write
\begin{align*}
&\max\{d(G(F(x_{2n},y_{2n}))),F(\alpha,\alpha')),d(G(F(y_{2n},x_{2n})),F(\alpha',\alpha))\}\\
&\leq \max\{d(G(G y_{2n}),S(\alpha')),d(G(G x_{2n}),S(\alpha))\}.
\end{align*}
Letting $n\rightarrow +\infty$, using the continuity of $G$, we obtain
\begin{equation*}
\max\{d(G(\alpha),F(\alpha,\alpha')),d(G(\alpha'),F(\alpha',\alpha))\}\leq \max\{d(G(\alpha),S(\alpha)),d(G(\alpha'),S(\alpha'))\}.
\end{equation*}
Looking at (\ref{eq14}), we deduce that
\begin{equation*}
\max\{d(G(\alpha),F(\alpha,\alpha')),d(G(\alpha'),F(\alpha',\alpha))\}=0.
\end{equation*}
Therefore,
\begin{equation*}
d(G(\alpha),F(\alpha,\alpha'))=0\quad\textrm{and}\quad d(G(\alpha'),F(\alpha',\alpha))=0.
\end{equation*}
Consequently
\begin{equation}\label{eq17}
G(\alpha)=F(\alpha,\alpha')\quad\textrm{and}\quad G(\alpha')=F(\alpha',\alpha).
\end{equation}
By the same way, we get
\begin{equation}\label{eq18}
S(\alpha)=F(\alpha,\alpha')\quad\textrm{and}\quad S(\alpha')=F(\alpha',\alpha).
\end{equation}
Finally, combining (\ref{eq14}), (\ref{eq17}) and (\ref{eq18}), we deduce that $(\alpha,\alpha')$ is a coupled coincidence point of $F$, $G$ and $S$.
\hfill $\blacksquare$

Now, we give a sufficient condition for the existence and the uniqueness of the coupled common fixed point. Notice that if $(X,\preceq)$ is a partially ordered set, we endow $X\times X$ with the following partial order relation:
$$\textrm{for }\,\,(x,y),\,(u,v)\,\in\,X\times X,\quad (x,y)\preceq (u,v)\Leftrightarrow x\preceq u\,\,\,\textrm{and}\,\,\,y\succeq  v.$$

\begin{Theorem}\label{T3}
In addition to the hypotheses of Theorem~{\ref{T1}} (resp. Theorem~{\ref{T2}}), suppose that for every $(x,\,y),\,(x^{*},\,y^{*})\,\in\,X\times X$ there exists a $(u,\,v)\,\in\,X\times X$ such that $(F(u,\,v),\,F(v,\,u))$ is comparable to $(F(x,\,y),\,F(y,\,x))$ and $(F(x^{*},\,y^{*}),\,F(y^{*},\,x^{*}))$. Then $F$, $G$ and $S$ have a unique coupled common fixed point, that is, there exist a unique $(x,\,y)\,\in\,X\times X$ such that
$$x=G(x)=F(x,\,y)=S(x)\quad\textrm{and}\quad y=G(y)=F(y,\,x)=S(y).$$
\end{Theorem}
{\bf Proof.}
We know, from Theorem~{\ref{T1}} (resp. Theorem~{\ref{T2}}), that exists a coupled coincidence point. We suppose that exist $(x,\,y)$ and $(x^{*},\,y^{*})$ two coupled coincidence points, that is, $G(x)=F(x,\,y)=S(x)$, $G(y)=F(y,\,x)=S(y)$, $G(x^{*})=F(x^{*},\,y^{*})=S(x^{*})$ and $G(y^{*})=F(y^{*},\,x^{*})=S(y^{*})$.\\
We claim that
\begin{equation}\label{eq19}
G(x)=G(x^{*})=S(x^{*})=S(x)\quad\textrm{and}\quad G(y)=G(y^{*})=S(y^{*})=S(y).
\end{equation}
If $(F(x,\,y),\,F(y,\,x))$ is comparable to $(F(x^{*},\,y^{*}),\,F(y^{*},\,x^{*}))$, it is easy to reach the result, then we suppose the general case.\\
By assumption there is $(u,\,v)\,\in\, X\times X$ such that $(F(u,\,v),\,F(v,\,u))$ is comparable to $(F(x,\,y)\,F(y,\,x))$ and $(F(x^{*},\,y^{*})\,F(y^{*},\,x^{*}))$. We distinguish two cases:\\
{\bf First case:} We assume that\\
$(F(x,\,y),\,F(y,\,x))\preceq (F(u,\,v),\,F(v,\,u))$ and $(F(x^{*},\,y^{*}),\,F(y^{*},\,x^{*}))\preceq (F(u,\,v),\,F(v,\,u))$.\\
Put $u_{0}=u$ and $v_{0}=v$ and we choose $u_{1}$ and $v_{1}$ such that $G(u_{0})\preceq S(u_{1})\preceq F(u_{0},\,v_{0})$, $G(v_{0})\succeq S(v_{1})\succeq F(v_{0},\,u_{0})$.\\
Similarly as in the proof of Theorem~{\ref{T1}}, we can construct sequences $\{u_{n}\}$ and $\{v_{n}\}$ in $X$ such that
\begin{equation*}
\begin{array}{ccc}
\left\{\begin{array}{ll}
G(u_{2n+2})=F(u_{2n},v_{2n})\\
G(v_{2n+2})=F(v_{2n},u_{2n})
\end{array}
\right.
\quad\textrm{and}\quad
\left\{\begin{array}{ll}
S(u_{2n+3})=F(u_{2n+1},v_{2n+1})\\
S(v_{2n+3})=F(v_{2n+1},u_{2n+1})
\end{array}
\right.
\quad\textrm{for all $n\geq 0$.}
\end{array}
\end{equation*}
Looking at the proof of Theorem~{\ref{T1}}, precisely at (\ref{eq2}), we see that $\{G(u_{2n})\}$ is a non-decreasing sequence, $G(u_{2n})\leq S(u_{2n+1})$, and $\{G(v_{2n})\}$ is a non-increasing sequence, $G(v_{2n})\succeq S(v_{2n+1})$.\\
Therefore, we have
\begin{eqnarray}\label{eq20}
G(x)=F(x,y)\leq F(u_{0},v_{0})=G(u_{2})\preceq G(u_{2n})\preceq S(u_{2n+1})\nonumber\\[\medskipamount]
\textrm{and}\hspace{20mm}\\[\medskipamount]
G(y)=F(y,x)\succeq F(v_{0},u_{0})=G(v_{2})\succeq G(v_{2n})\succeq S(v_{2n+1}).\nonumber
\end{eqnarray}
Similarly, we have
\begin{eqnarray}\label{eq21}
G(x^{*})=F(x^{*},y^{*})\preceq F(u_{0},v_{0})=G(u_{2})\preceq G(u_{2n})\preceq S(u_{2n+1})\nonumber\\[\medskipamount]
\textrm{and}\hspace{20mm}\\[\medskipamount]
G(y^{*})=F(y^{*},x^{*})\succeq F(v_{0},u_{0})=G(v_{2})\succeq G(v_{2n})\succeq S(v_{2n+1}).\nonumber
\end{eqnarray}
Using (\ref{eq20}) and the contractive condition, we write
\begin{eqnarray*}
\varphi(d(F(x,y),F(u_{2n+1},v_{2n+1})))&\leq \varphi(\max\{d(Gx,S u_{2n+1}),d(Gy,S v_{2n+1})\})\\[\medskipamount]
&-\phi(\max\{d(Gx,S u_{2n+1}),d(Gy,S v_{2n+1})\})
\end{eqnarray*}
and
\begin{eqnarray*}
\varphi(d(F(y,x),F(v_{2n+1},u_{2n+1})))&\leq \varphi(\max\{d(Gy,S v_{2n+1}),d(Gx,S u_{2n+1})\})\\[\medskipamount]
&-\phi(\max\{d(Gy,S v_{2n+1}),d(Gx,S u_{2n+1})\}).
\end{eqnarray*}
Therefore
\begin{eqnarray*}
\varphi(\max\{d(F(x,y),F(u_{2n+1},v_{2n+1})),d(F(y,x),F(v_{2n+1},u_{2n+1}))\})\\[\medskipamount]
\leq \varphi(\max\{d(G x,S u_{2n+1}),d(G y,S v_{2n+1})\})\\[\medskipamount]
-\phi(\max\{d(G x,S u_{2n+1}),d(G y,S v_{2n+1})\}).
\end{eqnarray*}
Therefore
\begin{eqnarray}\label{eq22}
\varphi(\max(d(G(x),S u_{2n+3}),d(Gy,S v_{2n+3})))&\leq \varphi(\max(d(G(x),S u_{2n+1}),d(Gy,S v_{2n+1})))\\[\medskipamount]
&-\phi(\max(d(G(x),S u_{2n+1}),d(Gy,S v_{2n+1}))).\nonumber
\end{eqnarray}
We see that
\begin{equation*}
\varphi(\max\{d(Gx,S u_{2n+3}),d(Gy,S v_{2n+3})\})\leq \varphi(\max\{d(Gx,S u_{2n+1}),d(Gy,S v_{2n+1})\}).
\end{equation*}
Using the non-decreasing property of $\varphi$, we get
\begin{equation*}
\max\{d(Gx,S u_{2n+3}),d(Gy,S v_{2n+3})\}\leq \max\{d(Gx,S u_{2n+1}),d(Gy,S v_{2n+1})\}.
\end{equation*}
This implies that $\max\{d(Gx,S u_{2n+1}),d(Gy,S v_{2n+1})\}$ is a non-increasing sequence.\\
Hence, there exists $r\geq 0$ such that
$$ \lim_{n\rightarrow +\infty}\max\{d(Gx,S u_{2n+1}),d(Gy,S v_{2n+1})\}=r.$$
Passing to limit in (\ref{eq22}) as $n\rightarrow +\infty$, we obtain
$$\varphi(r)\leq \varphi (r)-\phi(r),$$
which implies that $\phi(r)=0$ and then, since $\phi$ is an altering distance function, $r=0$.\\
We deduce that
\begin{equation}\label{eq23}
\lim_{n\rightarrow +\infty}\max\{d(Gx,S u_{2n+1}),d(Gy,S v_{2n+1})\}=0.
\end{equation}
Similarly, one can prove that
\begin{equation}\label{eq24}
\lim_{n\rightarrow +\infty}\max\{d(Gx^{*},S u_{2n+1}),d(Gy^{*},S v_{2n+1})\}=0.
\end{equation}
By the triangle inequality, (\ref{eq23}) and (\ref{eq24}),
\begin{eqnarray}
d(Gx,\,Gx^{*})\leq d(Gx,S u_{2n+1})+d(G(x^{*}),S u_{2n+1})\rightarrow 0\quad\textrm{as}\,\,n\rightarrow +\infty,\\[\medskipamount]
d(Gy,\,Gy^{*})\leq d(Gy,S v_{2n+1})+d(G(y^{*}),S v_{2n+1})\rightarrow 0\quad\textrm{as}\,\,n\rightarrow +\infty.
\end{eqnarray}
Hence
\begin{equation}\label{eqlak}
G(x)=G(x^{*})\,\,\,\textrm{and}\,\,\, G(y)=G(y^{*}).
\end{equation}
This prove the claim (\ref{eq19}) in this case.\\
{\bf Second case:} We assume that
$(F(x,\,y),\,F(y,\,x))\succeq (F(u,\,v),\,F(v,\,u))$ and $(F(x^{*},\,y^{*}),\,F(y^{*},\,x^{*}))\succeq (F(u,\,v),\,F(v,\,u))$.\\
Put $u_{0}=u$ and $v_{0}=v$ and we choose $u_{1}$ and $v_{1}$ such that $G(u_{0})\succeq S(u_{1})\succeq F(u_{0},\,v_{0})$, $G(v_{0})\preceq S(v_{1})\preceq F(v_{0},\,u_{0})$.\\
Similarly as in the proof of Theorem~{\ref{T1}}, we can construct sequences $\{u_{n}\}$ and $\{v_{n}\}$ in $X$ such that
\begin{equation*}
\begin{array}{ccc}
\left\{\begin{array}{ll}
G(u_{2n+2})=F(u_{2n},v_{2n})\\
G(v_{2n+2})=F(v_{2n},u_{2n})
\end{array}
\right.
\quad\textrm{and}\quad
\left\{\begin{array}{ll}
S(u_{2n+3})=F(u_{2n+1},v_{2n+1})\\
S(v_{2n+3})=F(v_{2n+1},u_{2n+1})
\end{array}
\right.
\quad\textrm{for all $n\geq 0$.}
\end{array}
\end{equation*}
Looking at the proof of Theorem~{\ref{T1}}, precisely at (\ref{eq2}), we see that $\{G(u_{2n})\}$ is a non-increasing sequence, $G(u_{2n})\succeq S(u_{2n+1})$, and $\{G(v_{2n})\}$ is a non-decreasing sequence, $G(v_{2n})\preceq S(v_{2n+1})$.\\
Therefore, we have
\begin{eqnarray*}
G(x)=F(x,y)\succeq F(u_{0},v_{0})=G(u_{2})\succeq G(u_{2n})\succeq S(u_{2n+1})\nonumber\\[\medskipamount]
\textrm{and}\hspace{20mm}\\[\medskipamount]
G(y)=F(y,x)\preceq F(v_{0},u_{0})=G(v_{2})\preceq G(v_{2n})\preceq S(v_{2n+1}).\nonumber
\end{eqnarray*}
Similarly, we have
\begin{eqnarray*}
G(x^{*})=F(x^{*},y^{*})\succeq F(u_{0},v_{0})=G(u_{2})\succeq G(u_{2n})\succeq S(u_{2n+1})\nonumber\\[\medskipamount]
\textrm{and}\hspace{20mm}\\[\medskipamount]
G(y^{*})=F(y^{*},x^{*})\preceq F(v_{0},u_{0})=G(v_{2})\preceq G(v_{2n})\preceq S(v_{2n+1}).\nonumber
\end{eqnarray*}
From this, we complete the proof identically  as in the first case and we obtain the claim (\ref{eq19}) in this case.
Since $G(x)=F(x,y)=S(x)$ and $G(y)=F(y,x)=S(y)$, by the commutativity of $F$, $G$ and $F$, $S$,  we have
\begin{eqnarray}\label{eq25}
\begin{array}{ccc}
\left\{
\begin{array}{lll}
G(G(x))=G(F(x,y))=F(Gx,Gy)\\
G(G(y))=G(F(y,x))=F(Gy,Gx)
\end{array}
\right.
\,\,\textrm{and}\,\,
\left\{
\begin{array}{lll}
S(S(x))=S(F(x,y))=F(S(x),S(y))\\
S(S(y))=S(F(y,x))=F(S(y),S(x)).
\end{array}
\right.
\end{array}
\end{eqnarray}
Set $G(x)=a=S(x)$, $G(y)=b=S(y)$. Then from (\ref{eq25}),
\begin{equation}\label{eq26}
G(a)=F(a,b)=S(a)\,\,\,\textrm{and}\,\,\,G(b)=F(b,a)=S(b).
\end{equation}
Thus $(a,b)$ is a coupled coincidence point. Then from (\ref{eq19}) with $x^{*}=a$ and $y^{*}=b$ it follows that $G(a)=G(x)=S(a)$ and $G(b)=G(y)=S(b)$.
Therefore
\begin{equation}\label{eq27}
G(a)=a=S(a)\quad\textrm{and}\quad G(b)=b=S(b).
\end{equation}
We deduce that $(a,b)$ is a coupled common fixed point. To prove the uniqueness, assume that $(c,d)$ is another coupled common fixed point. Then by (\ref{eq19}) and (\ref{eq27}) we have $c=G(c)=G(a)=a$ and $d=G(d)=G(b)=b$.
\hfill    $\blacksquare $

\begin{Remark}$\mbox{ }$\\
Taking $G=S=I_X$ (the identity mapping of $X$) in Theorem \ref{T1}, we obtain \cite[Theorem 2]{Harjani3}. \\
Taking $G=S=I_X$ in Theorem \ref{T2}, we obtain \cite[Theorem 3]{Harjani3}.
\end{Remark}

Taking $S=G$ in Theorem~\ref{T3},  we obtain the following result.
\begin{Corollary}\label{C1}
Let $(X,\,\preceq)$ be a partially ordered set and suppose that there exists a metric $d$ on $X$ such that $(X,d)$ is a complete metric space. Let $G: X\rightarrow X$ be two mappings and $F: X\times X\rightarrow X$ be a mapping with the mixed $G$-monotone property and satisfying
\begin{eqnarray*}
\varphi(d(F(x,y),F(u,v)))\leq \varphi(\max \{d(Gx,Gu),d(Gy,Gv)\})-\phi(\max\{Gx,Gu),d(Gy,Gv)\}),
\end{eqnarray*}
for all $x,\,y,\,u,\,v\,\in X$ with $G(x)\preceq G(u)$ or $G(x)\succeq G(u)$ and $G(y)\succeq G(v)$ or $G(y)\preceq G(v)$, where $\varphi$ and $\phi$ are altering distance functions.
 Assume that $F(X\times X)\subseteq G(X)$ and assume also the following hypotheses:
\begin{enumerate}
\item $G$ is continuous,
\item $F$ is continuous or $G$ is non-decreasing mapping and $X$ satisfies the following properties:
\begin{itemize}
\item if $(x_{n})$ is a non-decreasing sequences with $x_{n}\rightarrow x$ then $x_{n}\preceq x$ for each $n\in \mathbb{N}$,
\item if $(y_{n})$ is a non-increasing sequences with $y_{n}\rightarrow y$ then $y\preceq y_{n}$ for each $n\in \mathbb{N}$;\\
\end{itemize}
\item for every $(x,\,y),\,(x^{*},\,y^{*})\,\in\,X\times X$ there exists a $(u,\,v)\,\in\,X\times X$ such that $(F(u,\,v),\,F(v,\,u))$ is comparable to $(F(x,\,y),\,F(y,\,x))$ and $(F(x^{*},\,y^{*}),\,F(y^{*},\,x^{*}))$,
\item $F$ commutes with $G$.
\end{enumerate}
If there exist $x_{0},\,y_{0}\in X$  such that
\begin{equation*}
\left\{\begin{array}{ll}
G(x_{0})\preceq F(x_{0},y_{0})\\
G(y_{0}))\succeq F(y_{0},x_{0})
\end{array}
\right.
\end{equation*}
then there exists a unique $(x,\,y)\in X\times X$ such that
$$x=G(x)=F(x,y)\quad\textrm{and}\quad y=G(y)=F(y,x) ,$$
 that is, $G$ and $F$ have a unique coupled common fixed point.
\end{Corollary}

\section{Applications to periodic boundary value problems}

In this section, we study the existence and uniqueness of  solution to a periodic boundary value problem, as an application to the fixed point theorem given by Corollary~\ref{C1}. \\
Let $C([0,T],\mathbb{R})$ be the set  of all continuous functions $u\,:\,[0,T]\rightarrow \mathbb{R}$ and consider  a mapping $G\,:\,C([0,T],\mathbb{R})\rightarrow C([0,T],\mathbb{R})$.\\
Consider the periodic boundary value problem
\begin{align}
&u'=f(t,u)+h(t,u),\,\,\,t\in (0,T)\label{ap1}\\
&u(0)=u(T),
\label{ap2}
\end{align}
where $f$, $h$ are two continuous functions satisfying the following conditions:\\
There exist positive constants $\lambda_{1},\,\lambda_{2},\,\mu_{1}$ and $\mu_{2}$, such that for all $u,\,v\in (C([0,T],\mathbb{R})$, $ Gv(t)\leq Gu(t)$,
\begin{equation}
0\leq (f(t, u(t))+\lambda_{1} u(t))-(f(t, v(t))+\lambda_{1} v(t))\leq \mu_{1} \ln[(Gu(t)-Gv(t)
)^{2}+1]
\label{ap3}
\end{equation}
\begin{equation}
-\mu_{2}\ln[(Gu(t)-Gv(t))^{2}+1]\leq (h(t, u(t))+\lambda_{2} u(t))-(h(t, v(t))+\lambda_{2} v(t))\leq 0
\label{ap4}
\end{equation}
with
\begin{equation}\label{im}
\displaystyle\frac{2\max\{\mu_{1},\mu_{2}\}}{\lambda_{1}+\lambda_{2}}<1.
\end{equation}

We firstly study the existence of a solution of the following periodic system:
\begin{eqnarray}
u'+\lambda_{1}u-\lambda_{2}v=f(t,u)+h(t, v)+\lambda_{1}u-\lambda_{2}v\nonumber\\
v'+\lambda_{1}v-\lambda_{2}u=f(t, v)+h(t,u)+\lambda_{1}v-\lambda_{2}u,\label{ap5}
\end{eqnarray}
with the periodicity condition
\begin{equation}
u(0)=u(T)\quad\textrm{and}\quad v(0)=v(T).
\label{ap6}
\end{equation}
This problem is equivalent to the integral equations:
\begin{eqnarray*}
u(t)=\int_{0}^{T}k_{1}(t,s)[f(s,u)+h(s, v)+\lambda_{1}u-\lambda_{2}v]+\int_{0}^{T}k_{2}(t,s)[f(s, v)+h(s,u)+\lambda_{1}v-\lambda_{2}u]ds
\end{eqnarray*}
\begin{equation*}
v(t)=\int_{0}^{T}k_{1}(t,s)[f(s,v)+h(s, u)+\lambda_{1}v-\lambda_{2}u]+\int_{0}^{T}k_{2}(t,s)[f(s,u)+h(s,v)+\lambda_{1}u-\lambda_{2}v]ds
\end{equation*}
where
\begin{eqnarray*}
k_{1}(t,s)=\left\{\begin{array}{ll}\frac{1}{2}\left[\displaystyle\frac{e^{\sigma_{1}(t-s)}}{1-e^{\sigma_{1}T}}+\frac{e^{\sigma_{2}(t-s)}}{1-e^{\sigma_{2}T}}\right]
\quad &0\leq s< t\leq T \\[\bigskipamount]
\frac{1}{2}\left[\displaystyle\frac{e^{\sigma_{1}(t+T-s)}}{1-e^{\sigma_{1}T}}+\frac{e^{\sigma_{2}(t+T-s)}}{1-e^{\sigma_{2}T}}\right]
\quad &0\leq t< s\leq T
\end{array}
\right.
\end{eqnarray*}
\begin{eqnarray*}
k_{2}(t,s)=\left\{\begin{array}{ll}\frac{1}{2}\left[\displaystyle\frac{e^{\sigma_{2}(t-s)}}{1-e^{\sigma_{2}T}}+\frac{e^{\sigma_{1}(t-s)}}{1-e^{\sigma_{1}T}}\right]
\quad &0\leq s< t\leq T \\[\bigskipamount]
\frac{1}{2}\left[\displaystyle\frac{e^{\sigma_{2}(t+T-s)}}{1-e^{\sigma_{2}T}}+\frac{e^{\sigma_{1}(t+T-s)}}{1-e^{\sigma_{1}T}}\right]
\quad &0\leq t< s\leq T.
\end{array}
\right.
\end{eqnarray*}
Here, $\sigma_{1}=-(\lambda_{1}+\lambda_{2})$ and $\sigma_{2}=(\lambda_{2}-\lambda_{1})$.\\
From \cite[Lemma 3.2]{Bhaskar}, we have
\begin{equation}
k_{1}(t,s)\geq 0,\,\,\,0\leq t,s\leq T\quad\textrm{and}\quad k_{2}(t,s)\leq 0,\,\,\, 0\leq t,s\leq T.
\label{ap7}
\end{equation}
We assume that there exist $\alpha, \beta \in C([0,T])$ such that
\begin{eqnarray}
G(\alpha(t))\leq \int_{0}^{1}k_{1}(t,s)(f(s,\alpha(s))+h(s,\beta(s))+\lambda_{1}\alpha(s)-\lambda_{2}\beta(s))ds\nonumber\\[\medskipamount]
+\int_{0}^{1}k_{2}(t,s)(f(s,\beta(s))+h(s,\alpha(s))+\lambda_{1}\beta(s)-\lambda_{2}\alpha(s))ds\nonumber\\[\medskipamount]
\label{ap8}
\end{eqnarray}
and
\begin{eqnarray}
G(\beta(t))\geq \int_{0}^{1}k_{1}(t,s)(f(s,\beta(s))+h(s,\alpha(s))+\lambda_{1}\beta(s)-\lambda_{2}\alpha(s))ds\nonumber\\[\medskipamount]
+\int_{0}^{1}k_{2}(t,s)(f(s,\alpha(s))+h(s,\beta(s))+\lambda_{1}\alpha(s)-\lambda_{2}\beta(s))ds \label{ap9}.
\end{eqnarray}

We endow $X=C([0,T], \mathbb{R})$  with the metric $d(u,\,v)=\displaystyle\max_{t\in[0,T]}|u(t)-v(t)|$ for $u,\,v\in X$.\\
This space can be equipped with a partial order given by
$$x,\,y\in C([0,T]),\quad x\preceq y\Leftrightarrow x(t)\leq y(t),\quad\textrm{for any $t\in[0,T]$}.$$
In $X\times X$ we define the following partial order
$$(x,\,y),\,\,(u,\,v)\in X\times X,\quad (x,\,y)\preceq (u,\,v)\Leftrightarrow x\preceq u\,\quad \textrm{ and }\quad y\succeq v.$$
Since for any $x,\,y\in X$ we have that $\max(x,y)$ and $\min(x,y)\in X$, assumption~3 of Corollary~\ref{C1} is satisfied for $(X,\preceq)$.
Moreover in \cite{Nieto1} it is proved that $(X,\preceq)$ satisfies assumption~2 of Corollary~{\ref{C1}}.

Now, we shall prove the following result.
\begin{Theorem}
Suppose that $G: X\rightarrow X$ is a non-decreasing continuous mapping. Suppose also that
(\ref{ap3})-(\ref{im}) and  (\ref{ap8})-(\ref{ap9}) hold. Then (\ref{ap5})-(\ref{ap6}) has a unique solution. Therefore (\ref{ap1})-(\ref{ap2}) has also a unique solution.
\end{Theorem}
\noindent {\bf Proof.} We introduce the operator $F\,:\,X\times X\rightarrow X$ defined by
\begin{eqnarray*}
F(u,v)(t)&=&\int_{0}^{T}k_{1}(t,s)[f(s,u)+h(s, v)+\lambda_{1}u-\lambda_{2}v]\,\,ds\\
&&+\int_{0}^{T}k_{2}(t,s)[f(s, v)+h(s,u)+\lambda_{1}v-\lambda_{2}u]ds
\end{eqnarray*}
for all $u,v\in X$ and $t\in [0,T]$.

We claim that $F$ has the mixed $G$-monotone property.\\
In fact, for $Gx_{1}\leq Gx_{2}$ and $t\in[0,T]$, we have
\begin{eqnarray*}
F(x_{1},y)(t)-F(x_{2},y)(t)&=&\int_{0}^{T}k_{1}(t,s)(f(s,x_{1}(s))-f(s,x_{2})+\lambda_{1}(x_{1}(s)-x_{2}(s))ds\\
&&+\int_{0}^{T}k_{2}(t,s)(h(s,x_{1}(s))-h(s,x_{2})-\lambda_{2}(x_{1}-x_{2}))ds.
\end{eqnarray*}
From (\ref{ap3}), (\ref{ap4}) and (\ref{ap7}), for all $t\in [0,T]$, we have
$$
F(x_{1},y)(t)-F(x_{2},y)(t)\leq 0.
$$
This implies that
$$
F(x_{1},y)\preceq F(x_{2},y).
$$
Also, for $Gy_{1}\preceq Gy_{2}$ and $t\in[0,T]$, we have
\begin{eqnarray*}
F(x,y_{1})(t)-F(x,y_{1})(t)&=&\int_{0}^{T}k_{1}(t,s)(h(s,y_{1}(s))-h(s,y_{2})-\lambda_{2}(y_{1}(s)-y_{2}(s))ds\\
&& +\int_{0}^{T}k_{2}(t,s)(f(s,y_{1}(s))-f(s,y_{2})+\lambda_{1}(y_{1}-y_{2}))ds.
\end{eqnarray*}
Looking at (\ref{ap3}), (\ref{ap4}) and (\ref{ap7}), for all $t\in[0,T]$,  we have
$$
F(x,y_{1})(t)-F(x,y_{2})(t)\geq 0,
$$
that is,
$$
F(x,y_{1})\geq F(x,y_{2}).
$$
Thus, we proved that $F$ has the mixed $G$-monotone property.

For $G(x)\preceq G(u)$ and $G(y)\succeq G(v)$, we have $F(x,y)\succeq F(u,v)$ and
\begin{align*}
&d(F(x,y),F(u,v))=\max_{t\in[0,T]}|F(x,y)(t)-F(u,v)(t)|\\
&=\max_{t\in[0,T]}(F(x,y)(t)-F(u,v)(t))\\
&=\max_{t\in[0,T]}\int_{0}^{T}k_{1}(t,s)[(f(s,x(s))-f(s,u(s))+\lambda_{1}(x-u))-(h(s,v(s))-h(s,y(s))-\lambda_{2}(y-v))]ds\\
&-\int_{0}^{T}k_{2}(t,s)[(f(s,v(s))-f(s,y(s))+\lambda_{1}(v-y))-(h(s,u(s))-h(s,x(s))-\lambda_{2}(u-x))]ds.
\end{align*}
Using (\ref{ap3}) and (\ref{ap4}) we get
\begin{align*}
&d(F(x,y),F(u,v))\\
&\leq \max_{t\in[0,T]}\int_{0}^{T}k_{1}(t,s)\bigg(\mu_{1} \ln[(G x(s)-Gu(s))^{2}+1]+\mu_{2}\ln[(Gy(s)-Gv(s))^{2}+1]\bigg)ds\\
&+\int_{0}^{T}(-k_{2}(t,s))\bigg(\mu_{1} \ln[(G v(s)-Gy(s))^{2}+1]+\mu_{2} \ln[(Gx(s)-Gu(s))^{2}+1]\bigg)ds\\
&\leq \max(\mu_{1},\mu_{2})\max_{t\in[0,T]}\int_{0}^{T}(k_{1}(t,s)-k_{2}(t,s))\ln[(G x(s)-Gu(s))^{2}+1]ds\\
&+\int_{0}^{T}(k_{1}(t,s)-k_{2}(t,s))\ln[(G y(s)-Gv(s))^{2}+1]ds.
\end{align*}
An easy computation yields
\begin{align*}
&d(F(x,y),F(u,v))\\
&\leq\bigg(\max_{t\in[0,T]}\int_{0}^{T}(k_{1}(t,s)-k_{2}(t,s))ds\bigg)\max(\mu_{1},\mu_{2})\left(\ln[(d(G x,Gu))^{2}+1]+\ln[(d(G y,Gv))^{2}+1] \right)\\
&\leq 2 \bigg(\max_{t\in[0,T]}\int_{0}^{T}(k_{1}(t,s)-k_{2}(t,s))ds\bigg)\max(\mu_{1},\mu_{2})\ln[(\max(d(G x,Gu),d(Gy,Gv)))^{2}+1]\\
&\leq 2\max(\mu_{1},\mu_{2})\max_{t\in[0,T]}\left|\int_{0}^{t}\frac{e^{\sigma_{1}(t-s)}}{1-e^{\sigma_{1}T}}ds+\int_{t}^{T}\frac{e^{\sigma_{1}(t+T-s)}}{1-e^{\sigma_{1}T}}ds\right|\ln[(\max(d(G x,Gu),d(Gy,Gv)))^{2}+1].
\end{align*}
After integrating, we get
\begin{eqnarray*}
d(F(x,y),F(u,v))\leq \frac{2\max(\mu_{1},\mu_{2})}{\lambda_{1}+\lambda_{2}}\ln[(\max(d(G x,Gu),d(Gy,Gv)))^{2}+1].
\end{eqnarray*}
From (\ref{im}), we obtain
\begin{eqnarray*}
d(F(x,y),F(u,v))\leq \ln[(\max(d(G x,Gu),d(Gy,Gv)))^{2}+1]
\end{eqnarray*}
which implies that
\begin{eqnarray*}
(d(F(x,y),F(u,v)))^{2}\leq (\ln[(\max(d(G x,Gu),d(Gy,Gv)))^{2}+1])^{2}.
\end{eqnarray*}
Then,
\begin{align*}
&(d(F(x,y),F(u,v)))^{2}\leq (\max(d(G x,Gu),d(Gy,Gv)))^{2}\\
&-\left[(\max(d(G x,Gu),d(Gy,Gv)))^{2}-(\ln[(\max(d(G x,Gu),d(Gy,Gv)))^{2}+1])^{2}\right].
\end{align*}
Set $\varphi(t)=t^{2}$ and $\phi(t)=t^{2}-\ln(t^{2}+1)$. Clearly $\varphi$ and $\phi$ are altering distance functions and from the above inequality, we obtain
$$
\varphi(d(F(x,y),F(u,v)))\leq \varphi(\max\{d(G x,Gu),d(Gy,Gv)\})-\phi((\max\{d(G x,Gu),d(Gy,Gv)\}))
$$
for all $x,y,u,v\in X$ such that $G(x)\preceq G(u)$ and $G(y)\succeq G(v)$.\\
Now, let $\alpha, \beta \in X$ be the functions given by (\ref{ap8}) and (\ref{ap9}).Then,  we have
$$
G(\alpha)\preceq F(\alpha,\beta)\quad\textrm{ and }\quad F(\beta,\alpha)\preceq G(\beta).
$$

Thus, we proved that all the required hypotheses of Corollary~\ref{C1} are satisfied. Hence,
$G$ and $F$ have a unique coupled fixed point $(u,v)\in X\times X$, that is, $(u,v)$ is the unique solution of (\ref{ap5})-(\ref{ap6}).
\hfill $\blacksquare $

\vspace{0.5cm}

\noindent Habib Yazidi\\
UNIVERSIT\'E DE TUNIS, DEPARTMENT OF MATHEMATICS, TUNIS COLLEGE OF SCIENCES AND TECHNIQUES, 5 AVENUE TAHA HUSSEIN, BP, 59, BAB MANARA, TUNIS.\\
\textit{E-mail address:} habib.yazidi@gmail.com

\end{document}